\documentclass[]{amsart}
\usepackage{amsmath,amssymb, amscd}
\usepackage{pdfsync}
\usepackage{graphicx,subfigure}
\usepackage{color}
\usepackage[color]{showkeys}

\newcommand{\z}{\mathbf{z}}
\newcommand{\mr}{\mathbf{r}}

\newcommand{\R}{\mathbb{R}}

\newcommand{\Z}{\mathbb{Z}}
\newcommand{\B}{\mathbb{B}}
\def\pmtwo#1#2#3#4{\left( \begin{array}{cc}#1&#2\\#3&#4\end{array}\right)}
\def\dbar{\overline{\partial}}
\newcommand{\ee}{\mathrm{e}}
\newcommand{\ii}{\mathrm{i}}

\begin{document}
\title[Hybrid approach for DS II]{High precision numerical approach for the Davey-Stewartson 
II equation for Schwartz class initial data}

\author[C.~Klein]{Christian Klein}
\address[C.~Klein]
{Institut de Math\'ematiques de Bourgogne
9 avenue Alain Savary, BP 47870, 21078 Dijon Cedex}
\email{christian.klein@u-bourgogne.fr}

\author[K.~McLaughlin]{Ken McLaughlin}

\address[K.~McLaughlin]
{Department of Mathematics,
1874 Campus Delivery,
Fort Collins, CO 80523-1874}
\email{kenmcl@rams.colostate.edu}

\author[N.~Stoilov]{Nikola Stoilov}

\address[N.~Stoilov]{Institut de Math\'ematiques de Bourgogne, UMR 5584\\
                Universit\'e de Bourgogne-Franche-Comt\'e, 9 avenue Alain Savary, 21078 Dijon
                Cedex, France}
    \email{Nikola.Stoilov@u-bourgogne.fr}

\begin{abstract}
We present an efficient high-precision numerical approach for the 
Davey-Stewartson (DS) II equation, treating initial data from the Schwartz class of smooth, rapidly 
decreasing functions. As with previous approaches, the presented code uses 
discrete Fourier transforms for the spatial dependence and Driscoll's composite Runge-Kutta method for the time dependence. Since the DS equation is a nonlocal, nonlinear Schr\"odinger equation with a singular  symbol 
for the nonlocality, standard Fourier methods in practice 
only reach accuracies of the order of $10^{-6}$ or less for typical 
examples. This was previously demonstrated for 
the defocusing case by comparison with a numerical approach for DS via inverse scattering. By applying a regularization to the singular symbol, originally developed for D-bar problems, the presented code is shown to 
reach machine precision. The code can treat integrable and 
non-integrable DS II equations.  Moreover, it has the same 
numerical complexity as existing codes for DS II.  Several examples 
for the integrable defocusing DS II equation are discussed as test cases. 

In an appendix by C.~Kalla, a doubly periodic solution to 
the defocusing  DS II equation is presented, providing a test for direct
DS codes based on Fourier methods.

\end{abstract}

\date{\today}

\subjclass[2000]{}

\thanks{This work was partially supported  by 
the ANR-FWF project ANuI - ANR-17-CE40-0035, the isite BFC project 
NAANoD, the ANR-17-EURE-0002 EIPHI and by the 
European Union Horizon 2020 research and innovation program under the 
Marie Sklodowska-Curie RISE 2017 grant agreement no. 778010 IPaDEGAN.
K.M. was supported in part by the National 
Science Foundation under grant DMS-1733967. 
}
\maketitle

\section{Introduction} 

The general Davey-Stewartson (DS) II system is
\begin{equation}
\label{DSgen}
\begin{array}{ccc}
i
\Psi_{t}+ \Psi_{xx}-\Psi_{yy}+2\rho\left(\beta \Phi+\left| \Psi \right|^{2}\right)\Psi & = & 0,
\\
\Phi_{xx}+\Phi_{yy}+2\left| \Psi \right|_{xx}^{2} & = & 0,
\end{array}
\end{equation}
where $\beta$ is a positive constant, indices denote partial derivatives,    $\rho$ takes values $\pm1$, and where $\Phi$ denotes a mean field. 
The systems (\ref{DSgen}) are of considerable importance in 
applications since they are a simplification of the 
Benney-Roskes \cite{BR}, or Zakharov-Rubenchik \cite{ZR} systems being 
`universal' models for the description of the interaction of short and 
long waves. They first appeared in the context of water waves 
\cite{DS,DR,AS,Lan} in the so-called modulational (Schr\"odinger) regime, 
i.e., in the study of the modulation of plane waves.
In \cite{CL1,CL2} is was shown, via  
diffractive geometric optics, that DS systems provide good 
approximate solutions to general quadratic hyperbolic systems. 
Furthermore, the DS systems appear  
in numerous  physical contexts as ferromagnetism \cite{Leb}, plasma 
physics \cite{MRZ}, and nonlinear optics \cite{NM}.      
The Davey-Stewartson systems can also be viewed as the 
two-dimensional version of the Zakharov-Schulman systems (see 
\cite{ZS1,ZS2,GS}.
For more details on DS and its applications the reader is referred to 
\cite{KSDS,KSint} where an abundance of references can be found.

Because of their importance in applications, many numerical 
approaches for nonlinear Schr\"odinger (NLS) equations have been 
developed. For Schwartz  class or periodic functions, Fourier 
spectral methods generally are the most efficient, and there is also 
a multitude of adapted time integration schemes, see \cite{etna,BKutz} for  a comparison and references. The DS system can be seen as a nonlocal 
NLS equation  where the nonlocality corresponds to a singular Fourier 
symbol, see (\ref{DSIIbis}) below. Several numerical approaches have been presented 
for DS II systems along these lines, see \cite{WW,BMS,KR,KS}. In \cite{KM} it was shown that the standard 
treatment of this singular symbol leads for typical examples to a 
numerical error of the order of at least $10^{-6}$. The goal of 
this paper is to present an efficient numerical approach for the DS II 
systems 
so that machine precision can be reached, with the same computational complexity as existing codes. 

If the second equation in (\ref{DSgen}) is solved by formally 
inverting the Laplace operator $\Delta=\partial_{xx}+\partial_{yy} $
with some boundary conditions at infinity, the DS II equation can be 
written as a nonlocal nonlinear Schr\"odinger equation (NLS),  
\begin{equation}\label{DSIIbis}
   i \epsilon\partial_{t}\psi+\epsilon^{2}\partial_{xx}\psi-\epsilon^{2}\partial_{yy}\psi+2\rho\Delta^{-1}\lbrack
    \left(\partial_{yy}+(1-2\beta)\partial_{xx}\right)\left|\psi\right|^{2}\rbrack\psi  =  0,
    \end{equation}
 which involves the zeroth order  nonlocal operator
$$\Delta^{-1}\lbrack
    \left(\partial_{yy}+(1-2\beta)\partial_{xx}\right)\rbrack.$$
We put $\xi=\xi_{1}+i\xi_{2}$ where $\xi_{1}$ and $\xi_{2}$ are 
the dual Fourier variables to $x$ and $y$ respectively. The Fourier 
transform of a function $\Phi$ as defined in (\ref{F}) is denoted by 
$\hat{\Phi}=\mathcal{F}\Phi$. In Fourier space, equation 
(\ref{DSIIbis}) takes the form
\begin{equation}
    i
    \hat{\Psi}_{t}+\frac{1}{2}(\xi^{2}+\bar{\xi}^{2})\hat{\Psi}-2\rho 
    \mathcal{F}\left[\mathcal{F}^{-1}\left(\beta\left\{\frac{\xi}{2\bar{\xi}}
    +\frac{\bar{\xi}}{2\xi}\right\}
    \mathcal{F}|\Psi|^{2}\right)\Psi
    +(\beta-1)|\Psi|^{2}\Psi\right]=0
    \label{DSfourier}.
\end{equation}
Note that DS II reduces to the standard cubic NLS equation in 1d for 
$y$-independent initial data. We remark that here, and throughout this paper, a function depending on a single complex variable  is not necessarily
holomorphic in this variable.  

The problematic term for a Fourier spectral method is in the 
nonlocality in (\ref{DSfourier}), i.e., with $\xi=\chi\exp(i\psi)$, 
$\chi,\psi\in\mathbb{R}$, 
\begin{equation}
    \frac{\xi}{2\bar{\xi}}+\frac{\bar{\xi}}{2\xi} = \cos(2\psi)
    \label{quo}.
\end{equation}
This means that the term does not have a well defined limit for $\chi\to 0$. 
Previous codes, for example discussed in \cite{KR} simply fix the value there which leads to an integrand with a discontinuity for $\chi=0$.  We will henceforth call these methods \textit{classical}.  For these classical methods, the discontinuity implies that the inverse discrete Fourier transform will not show  \emph{spectral convergence}, that is, an exponential 
decay of the numerical error with the resolution. It 
turns out that the precise value between $-1$ and $1$ chosen for 
$\cos(2\psi)$ will not affect the numerical accuracy. Therefore we 
chose it to be 0 in previous codes.

%
%
%
%
%
%
%
%
A problem in developing numerical approaches for the DS II equation 
for Schwartz class or periodic data has been the absence of suitable 
exact test solutions. Although the equation is completely integrable 
for $\beta=1$
and though there is consequently a wealth of exact solutions to 
choose from, none of these previously known explicit solutions 
appears suitable in this context: the solitons of the 1d NLS equation 
are solutions of the DS II equation called line solitons. But they 
either do not test the 2d aspect of DS II, or they hit the 
computational domain at an angle which leads to Gibbs' phenomena if 
Fourier methods are applied. There is also a 2d soliton, the 
\emph{lump}, which, however, shows algebraic decay towards infinity in all spatial 
directions. Such functions, in contrast to Schwartz class functions, 
in practice cannot be continued periodically as being 
smooth within the finite numerical precision. There is, however, a large class 
of quasi-periodic solutions to the integrable DS II equation (with $\beta = 1$), given in terms of multi-dimensional 
theta functions on general compact Riemann surfaces. They were 
constructed 
by Malanyuk \cite{Mal} and have been re-derived by Kalla via Fay's 
trisecant identity and discussed in \cite{Kal}. A numerical approach 
to evaluate these explicit solutions is presented in \cite{KK}. In 
general, these solutions are not exactly periodic, and certainly not 
in $x$ and $y$ (for example, they can be periodic on a 
non-rectangular lattice in the $(x,y)$ plane). This can lead (again) to problems with a Fourier approach.

\begin{figure}[htb!]
   \includegraphics[width=0.7\textwidth]{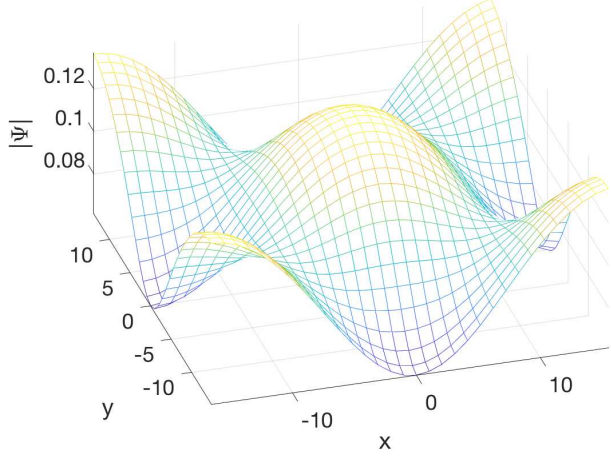}
\caption{Doubly periodic solution to the defocusing DS II equation 
discussed in the appendix for $t=0$.  }  
\label{figexact}
\end{figure}

In a previous work \cite{KR}, which compares  time integration schemes 
for DS II, the approach is to use a reference solution obtained 
from various such schemes with very high resolution. Whereas this 
works well to test the convergence of the time integration methods, 
this obviously does not test the potential mistakes, nor identify potential sources of errors, in the spatial 
discretization, since the latter is always the same  and the accuracy is 
unknown. In the present paper, 
we address this problem in two ways: first in Section 1 of the 
Appendix, written by Kalla, we present a family of solutions both 
periodic in $x$ and in $y$. As an example we consider a genus 2 
solution which is a travelling wave, see Fig.~\ref{figexact}. Such a solution was first given for the 
Kadomtsev-Petviashvili equation, a generalization to 2d of the Korteweg-de 
Vries equation with solutions similar  to DS II, in \cite{DFS}. This 
solution provides a useful test as shown in \cite{Dri,KR}. 
Secondly, we presented in \cite{KM,KMS} a high-precision numerical 
approach for the scattering problem for DS II which allows one to generate 
reference solutions in the defocusing case with essentially machine 
precision. Such high precision will be useful to study the stability 
(or the lack thereof) of certain exact solutions as solitons, and a 
potential blow-up in this context. 

Both types of tests are applied in this paper to show they 
are very useful if high precision approaches are looked for. 
Note that DS II equations are in general not integrable. However, for 
the tests of the codes we concentrate on the integrable DS II since 
we want to compare with a reference solution obtained via inverse 
scattering techniques as in \cite{KM}. The numerical techniques 
presented in this paper are directly applicable to the non-integrable 
case $\beta\neq1$.

In \cite{KM,KMS} we have shown how to regularize terms of the type 
(\ref{quo}) arising 
in the context of D-bar equations with a 
hybrid approach, i.e., by regularizing the term via analytically 
reducing them to terms regular within a finite 
numerical precision and then doing a numerical computation. Note that the terms treated in this way in 
\cite{KM,KMS} are less regular (actually unbounded) than the ones 
considered here. Therefore the regularization approach for DS II  is less 
computationally expensive than in the D-bar case. We show that this way we can reach machine precision even in cases for which this was not possible with previous codes.

This paper is organized as follows: in section \ref{Sec:HybNum} we give a brief 
summary of the numerical approaches for DS II via Fourier techniques 
and the regularization approaches of \cite{KM,KMS} applied to this 
case. In section \ref{Sec:Examples} we consider examples from the Schwartz class with 
and without radial symmetry which are compared to the numerical solution for the defocusing DS II for the same initial data constructed via the scattering approach. It is shown that machine precision can be essentially reached. In the first section of the appendix C. Kalla presents an exact solution 
to DS II periodic both in $x$ and $y$ on a genus 2 surface, and  in the second section of the appendix we summarize the scattering and inverse scattering theory of the DS II equation.

\section{Hybrid numerical approach for DS II}
\label{Sec:HybNum}
In this section we briefly describe the numerical approach for the DS 
II equation.

The Fourier transform of a function $\Phi$ is defined as 
\begin{align}
    \hat{\Phi} & = \mathcal{F}\Phi:=\frac{1}{2\pi}\int_{\mathbb{R}^{2}}^{}\Phi 
    e^{-i(\xi \bar{z}+\bar{\xi}z)/2}\Phi dx dy,
    \label{F}
    \\
    \Phi & =\mathcal{F}^{-1}\Phi=\frac{1}{2\pi}\int_{\mathbb{R}^{2}}^{} e^{i(\xi 
    \bar{z}+\bar{\xi}z)/2}\hat{\Phi} d\xi_{1}d\xi_{2}.
    \nonumber
\end{align}

The classical numerical approaches to DS II equations, see for instance 
\cite{WW,BMS,KR}, via Fourier spectral methods are based on an 
approximation of the Fourier transform (\ref{F}) by a discrete 
Fourier transform, i.e., by essentially a truncated 
Fourier series. For the latter an efficient algorithm exists known as 
the fast Fourier transform (FFT). This leads for (\ref{DSfourier}) to a finite 
dimensional system of ordinary differential equations (ODEs) for 
$\mathcal{F}\Psi$ (in an abuse of notation, we use the same symbol 
for the continuous and the discrete Fourier transform). The ODE system 
has the form
\begin{equation}
    \hat{\Psi}_{t}=\mathbf{L}\hat{\Psi}+\mathbf{N}[\hat{\Psi}]
    \label{sys},
\end{equation}
where $\mathbf{L}$ and $\mathbf{N}$ denote linear and nonlinear 
operators, respectively:
\begin{eqnarray}
&& \mathbf{L}\hat{\Psi} = - \frac{1}{2}(\xi^{2}+\bar{\xi}^{2})\hat{\Psi} \\
&&\mathbf{N}[\hat{\Psi}] = -2\rho 
    \mathcal{F}\left[\mathcal{F}^{-1}\left(\beta\left\{\frac{\xi}{2\bar{\xi}}
    +\frac{\bar{\xi}}{2\xi}\right\}
    \mathcal{F}|\Psi|^{2}\right)\Psi
    +(\beta-1)|\Psi|^{2}\Psi\right]
    \label{sys2}
\end{eqnarray}

 Since the operator $\mathbf{L}$ is 
proportional to $|\xi|^{2}$ whereas $\mathbf{N}$ is of order $|\xi|^{0}$, 
this is an example of a \emph{stiff} system of equations, in particular if 
large values of $|\xi|$ are needed as in the case of rapid 
oscillations called \emph{dispersive shock waves} as in \cite{KR2} 
or in the case of a \emph{blow-up}, a diverging $L^{\infty}$ norm of 
the solution as studied in \cite{KMR,KSDS,KS}. Loosely speaking the word 
\emph{stiff} refers to ODEs for which explicit time integration 
schemes are not efficient for stability reasons which means that the 
condition of numerical stability requires considerably smaller time 
steps than imposed by the desired accuracy. For equations of the form 
(\ref{sys}) with stiffness in the linear term, many efficient 
integrations schemes exist in the literature. Several such methods 
have been compared for DS II in \cite{KR} where the reader can find an abundance of references. The best method identified there for both the focusing 
and the defocusing case is Driscoll's composite Runge-Kutta (RK) method 
\cite{Dri}. It uses a stiffly stable third order RK method 
for the high wavenumbers in $\mathbf{L}$ and the standard explicit 
fourth order RK method for the remaining terms. We apply this method 
in the following.

The task is thus to compute two singular integrals of the form 
$\mathcal{F}^{-1}(S/\xi)$ and $\mathcal{F}^{-1}(S/\bar{\xi})$. In 
\cite{KM,KMS} we have shown that (we have $\partial = 
(\partial_{x}-i\partial_{y})/2$)
\begin{equation}
    \mathcal{F}^{-1}\left(\frac{\bar{\xi}^{n}e^{-|\xi|^{2}}}{\xi}\right)=
    (-2i\partial)^{n}\frac{i}{z}\left(1-e^{-\frac{|z|^{2}}{4}}
    \right),\quad n=0,1,\ldots
   \label{expxi2},
\end{equation}
and 
\begin{equation*}
    \mathcal{F}^{-1}\left(\frac{ \xi^{n}e^{-|\xi|^{2}}}{\bar{\xi}}\right)=
    (-2i\bar{\partial})^{n}\frac{i}{\bar{z}}\left(1-e^{-\frac{|z|^{2}}{4}}
    \right),\quad n=0,1,\ldots
\end{equation*}

The idea to compute a singular integral of the form 
$\mathcal{F}^{-1}\left(S(\xi)/\xi\right)$ is to subtract the first 
$M+1$ terms of a Taylor series of $S$  in  $\bar{\xi}$ where $M$ is chosen 
such that the residual is smooth up to the finite numerical 
precision, i.e., 
\begin{equation}
    \mathcal{F}^{-1}\left(\frac{S(\xi)}{\xi}\right)
    =     
    \mathcal{F}^{-1}\left(\frac{S(\xi)-e^{- 
    |\xi|^{2}}\sum_{n=0}^{M}\partial_{\bar{\xi}}^{n}
    S(0)\bar{\xi}^{n}/n!}{\xi}\right)+ e^{- |\xi|^{2}} \sum_{n=0}^{M}\partial_{\bar{\xi}}^{n}
    S(0) (-2i\partial)^{n}\frac{i}{z}\left(1-e^{-\frac{|z|^{2}}{4}}
    \right).
    \label{xireg1}
\end{equation}
The first term on the right hand side of (\ref{xireg1}) can now be 
efficiently computed numerically via FFT, the 
second is computed analytically. 
The derivatives of $S$ are also computed via Fourier techniques,
\begin{align}
    \partial_{\xi}^{n}
    S(\xi)&=\mathcal{F}\left[(i\bar{z}/2)^{n}\mathcal{F}^{-1}S\right]
    \nonumber\\
    \partial_{\bar{\xi}}^{n}
    S(\xi)&=\mathcal{F}\left[(-iz/2)^{n}\mathcal{F}^{-1}S\right]
    \label{xireg2}
\end{align}
Note that the derivatives are only needed for $\xi=0$. To compute the 
first term on the right-hand-side of (\ref{xireg1}) at $\xi=0$, one 
would  also need $\partial_{\xi}S(0)$ (according to L'H\^opital's rule). 

Integrals of the form  
$\mathcal{F}^{-1}\left(S(\xi)/\bar{\xi}\right)$ are computed 
analogously. As already mentioned, for DS II the situation is less 
singular than the cases discussed in \cite{KM,KMS}. Here we only need 
to compute the singular term in (\ref{DSfourier}) or, equivalently, the singular part of the nonlinear term appearing in (\ref{sys})-(\ref{sys2}). This leads to 
\begin{equation}
    \begin{split}
   & \mathcal{F}^{-1}\left(\left[\frac{\bar{\xi}}{2\xi}+
    \frac{\xi}{2\bar{\xi}}\right]S(\xi)\right)\\
    &=  \mathcal{F}^{-1}\left(\left[\frac{\bar{\xi}}{2\xi}+
    \frac{\xi}{2\bar{\xi}}\right]S(\xi)-\frac{1}{2}e^{- 
    |\xi|^{2}}\sum_{n=0}^{M}\left[\partial_{\bar{\xi}}^{n}
    S(0)(\bar{\xi}^{n+1}/\xi)+\partial_{\xi}^{n}
    S(0)(\xi^{n+1}/\bar{\xi})\right]/n!\right)\\
    &+ \frac{1}{2}e^{- |\xi|^{2}} \sum_{n=0}^{M}\left[\partial_{\bar{\xi}}^{n}
    S(0) (-2i\partial)^{n+1}\frac{i}{z}\left(1-e^{-\frac{|z|^{2}}{4}}
    \right)+\partial_{\xi}^{n}
    S(0) (2i\bar{\partial})^{n+1}\frac{-i}{\bar{z}}\left(1-e^{-\frac{|z|^{2}}{4}}
    \right)\right].
     \end{split}
    \label{DSreg}
\end{equation}

A brief summary of the algorithm to solve the DS II equation (\ref{DSgen}) is as follows.  
\begin{enumerate}
\item Introduce a discrete Fourier grid on the spatial domain $ x \in L_{x}[-\pi, \pi ], \ \ y \in L_{y}[-\pi, \pi] $.
\item Compute and store the functions 
\begin{eqnarray*}
W_{n} = \mathcal{F}^{-1} \left( \frac{ \overline{\xi}^{n}}{\xi} e^{- | \xi |^{2} } \right) \ .
\end{eqnarray*}
\item In the implementation of Driscoll's Runge-Kutta method, compute the nonlinear term $\mathbf{N}\left(\hat{\Psi} \right)$ in (\ref{sys})-(\ref{sys2}) as follows:
\begin{enumerate}
\item Compute the relevant Taylor coefficients of $ S( \xi ) = 
\mathcal{F}( | \Psi |^{2}) $ directly from $|\Psi|^{2}$ via 
(\ref{xireg2}) (i.e. without computing any Fourier transforms):
\begin{eqnarray*}
   \partial_{\bar{\xi}}^{n}
 \mathcal{F}( |\Psi|^{2})=\int_{\mathbb{R}^{2}} \left[(-iz/2)^{n} | 
 \Psi |^{2} \right] dx dy \ ,
    \label{xireg2}
\end{eqnarray*}
where the integral is simply approximated by a sum over all indices 
after the spatial discretization 
(note that the application of the trapezoidal rule in this context is 
a spectral method since it corresponds simply to the computation of 
the Fourier coefficient with indices 0). 

\item Compute $S(\xi) = \mathcal{F}\left(|\Psi|^{2}\right)$, then compute 
\begin{eqnarray*}&&
 \beta  \mathcal{F}^{-1}\left(\left[\frac{\bar{\xi}}{2\xi}+
    \frac{\xi}{2\bar{\xi}}\right]S(\xi)-\frac{1}{2}e^{- 
    |\xi|^{2}}\sum_{n=0}^{M}\left[\partial_{\bar{\xi}}^{n}
    S(0)(\bar{\xi}^{n+1}/\xi)+\partial_{\xi}^{n}
    S(0)(\xi^{n+1}/\bar{\xi})\right]/n!\right) ,
\end{eqnarray*}
using (a) above.
\item Now add to this quantity the exact term appearing in the last line of (\ref{DSreg}), and also the term $(\beta-1)|\Psi|^{2}$.
\item Finally, the nonlinear term $\mathbf{N}\left(\hat{\Psi}\right)$ is obtained by multiplying the result by $\Psi$, and computing the Fourier transform of that product.
\end{enumerate}
\end{enumerate}

The reader will note that the functions $W_{n}$ in (2) need only be 
computed once, they are the same for all time steps and for all 
stages of the RK scheme. The computation 
of the nonlinear term $\mathbf{N}\left( \hat{\Psi} \right)$ using the 
regularization described in (a) through (d) above then {\it requires 
the same number of Fourier transforms as the classical approach 
without regularization}.  Therefore, the numerical complexity of this algorithm is the same as the standard algorithm.

\section{Examples}
\label{Sec:Examples}
In this section we compare DS codes with and without the 
regularization of the nonlocal term of the previous section for 
various examples in the defocussing integrable case $\beta=1$. We first consider the exact solution, and then more 
general initial data with a numerical scattering approach providing 
the reference solution.

\subsection{Exact solution}
\begin{figure}[htb!]
   \includegraphics[width=0.7\textwidth]{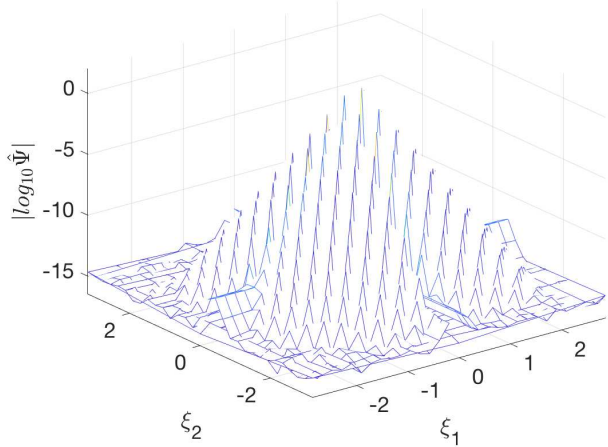}
\caption{Modulus of the Fourier coefficients of the doubly 
periodic solution to the defocusing DS II equation 
discussed in the appendix for $t=0$.  }  
\label{fourierexact}
\end{figure}

As a first test we take as initial data with $t=0$ the exact periodic solution 
of Fig.~\ref{figexact}. We use 
$N_{x}=N_{y}=2^{5}$ Fourier modes in $x$ and $y$ direction. In 
Fig.~\ref{fourierexact} we show the Fourier coefficients for this 
solution as computed with an FFT. It can be seen that the Fourier 
coefficients decrease in all directions to $10^{-15}$, the level of 
the rounding errors. This shows that the solution is numerically 
well resolved. 

As similar exact solutions to KdV and KP, this solution is numerically 
resolved with a comparatively small number of Fourier modes. As 
discussed for instance in \cite{KR}, this implies that the system 
(\ref{sys}) is not stiff, and that this example is thus not 
challenging for the codes.
We use $N_{t}=1000$ 
time steps for $t\leq10$.  The $L^{2}$ norm of the DS II solution, 
which is a conserved quantity for DS II, is not automatically 
conserved by the code and thus a test of the time resolution. In this 
example, it is conserved to the order of $10^{-16}$ during the entire 
computation. The difference to the exact solution at the final time 
is of the order of $10^{-14}$, the estimated accuracy for the exact 
solution.

This shows that the time evolution code reproduces the solution with 
the accuracy with which it is known. Note that this is both the case 
for the codes with and without the regularization of the previous 
section. Thus this exact solution does not test this aspect of the 
codes. In the next subsection  we thus present examples that require higher resolution in space.

\subsection{Scattering approach}

The defocussing DSII equation (\ref{DSgen}) with $\beta=1$ was shown to be integrable in  \cite{ASBook}.  Scattering and inverse scattering theory were developed in a number of works, including \cite{APP}, and \cite{Sung}. Asymptotic behavior of solutions was considered in \cite{Sung}, and more recently in \cite{Perry2012} and separately \cite{NRT}, asymptotic results were obtained with initial data in natural function spaces.  In Section 2 of the Appendix, we provide a brief summary of the scattering and inverse scattering theory that provides a solution procedure for the integrable case.

In \cite{KM} a spectral method is developed to solve the D-bar 
equations that appear in the scattering and inverse scattering theory 
associated to the integrable defocusing DS II equation.  The method 
is shown to have spectral convergence in the number of Fourier 
modes, for Schwartz class functions.  The approach is based on the 
integral equation in the Fourier domain associated to the D-bar system (\ref{eq:specprob2}), which contains a singular integral operator. The singular part appearing in the singular integral operator  is subtracted off using polynomials and exponentials, and then the offending terms are handled via exact calculations, just as we are implementing in (\ref{DSreg}) above.

We consider the computational method in \cite{KM} as a source of a {\it completely independent collection of reference solutions to the DS II equation}.  Here we use two examples to study the accuracy of the new method we present in Section \ref{Sec:HybNum}.

\begin{figure}[htb!]
\includegraphics[width=0.45\textwidth]{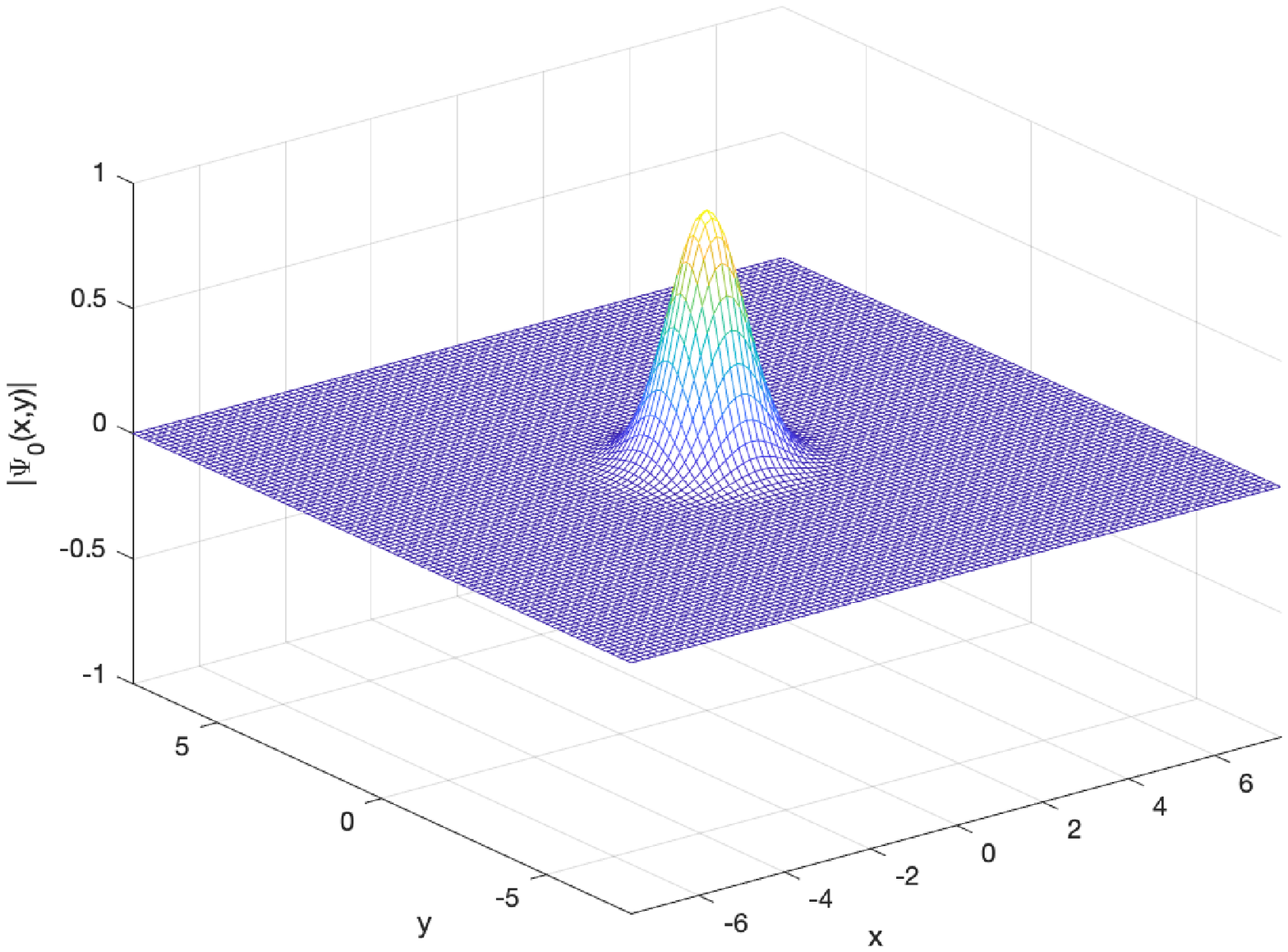}
\includegraphics[width=0.45\textwidth]{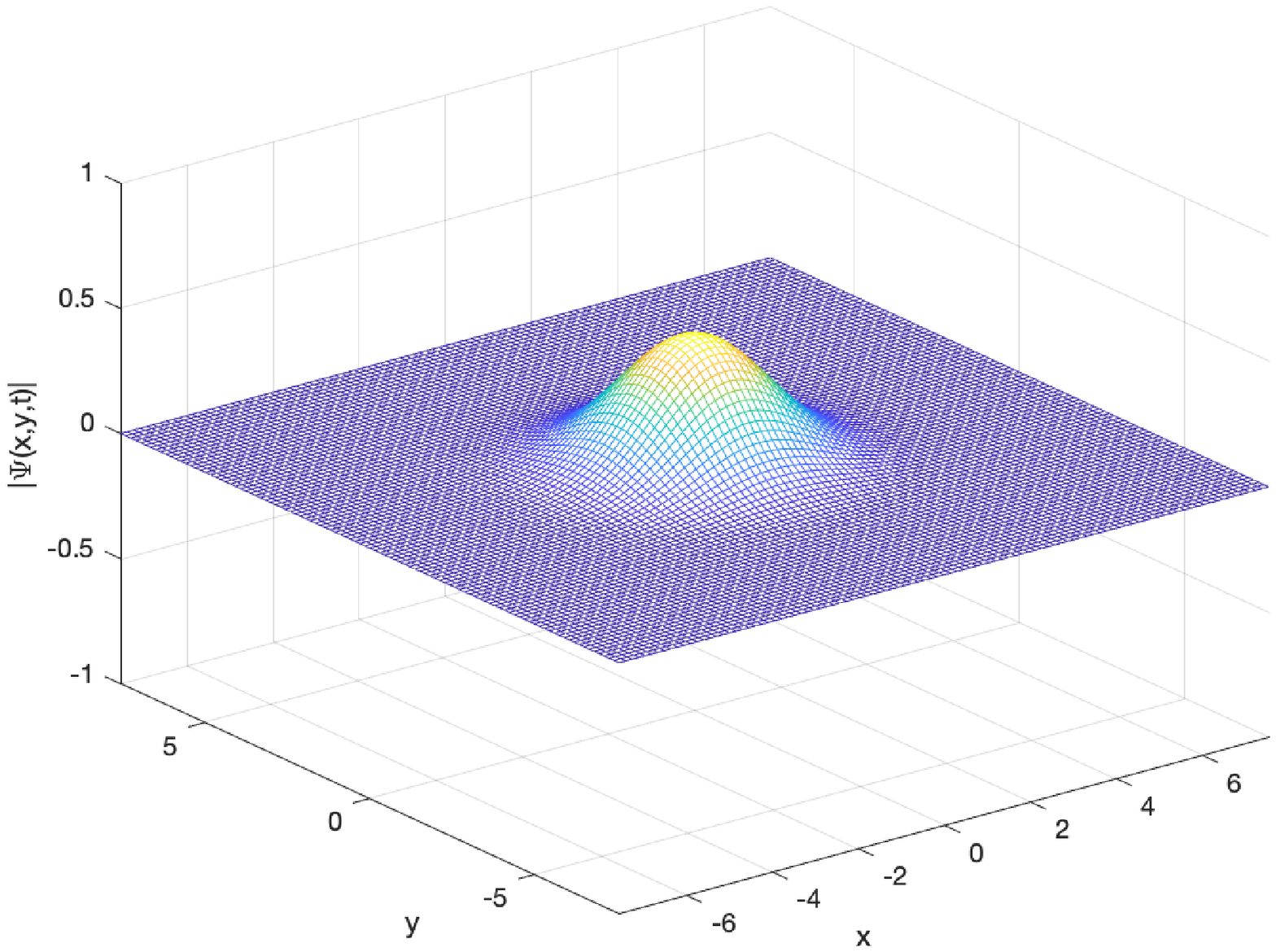}  \includegraphics[width=0.45\textwidth]{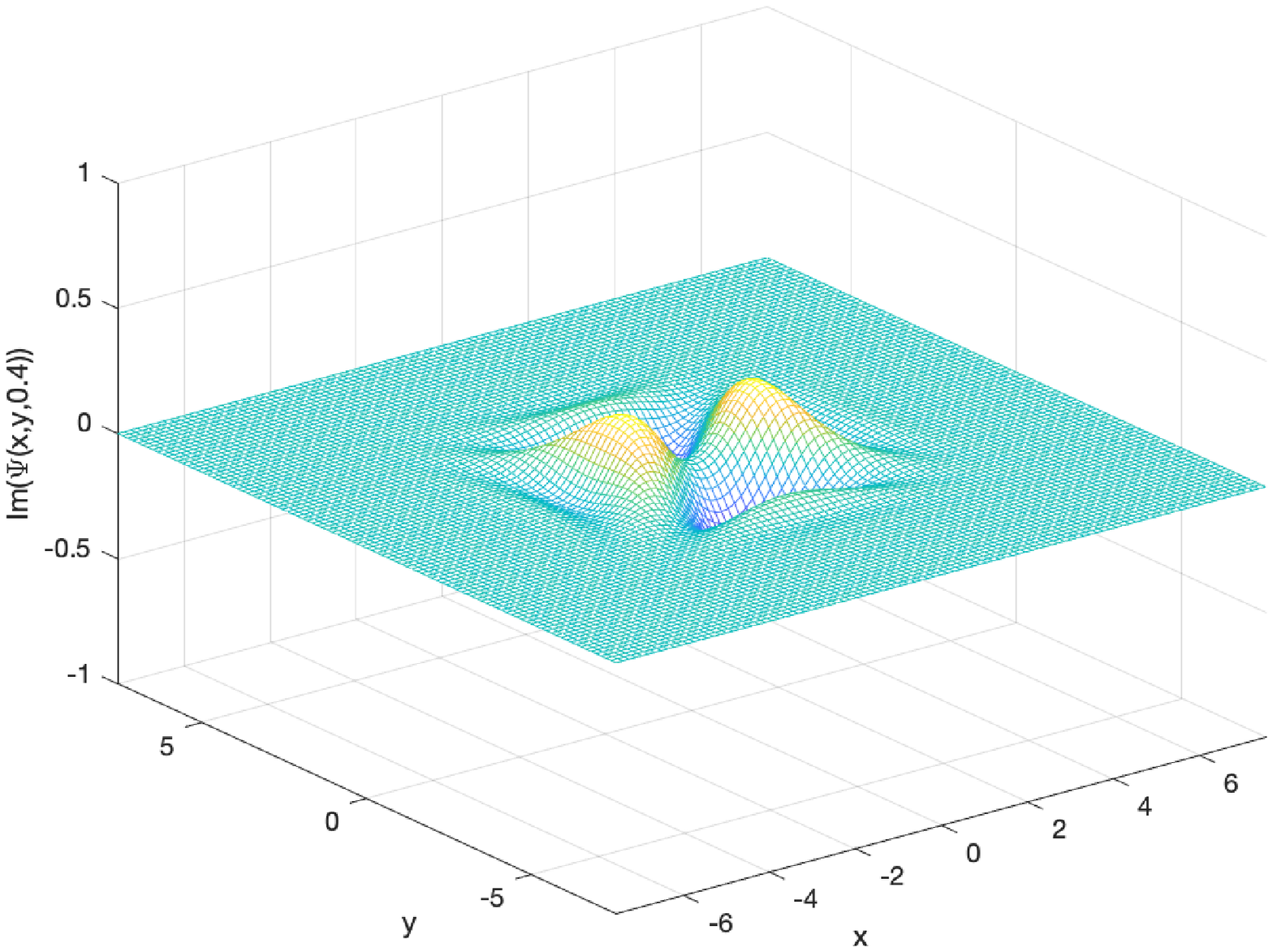}\includegraphics[width=0.45\textwidth]{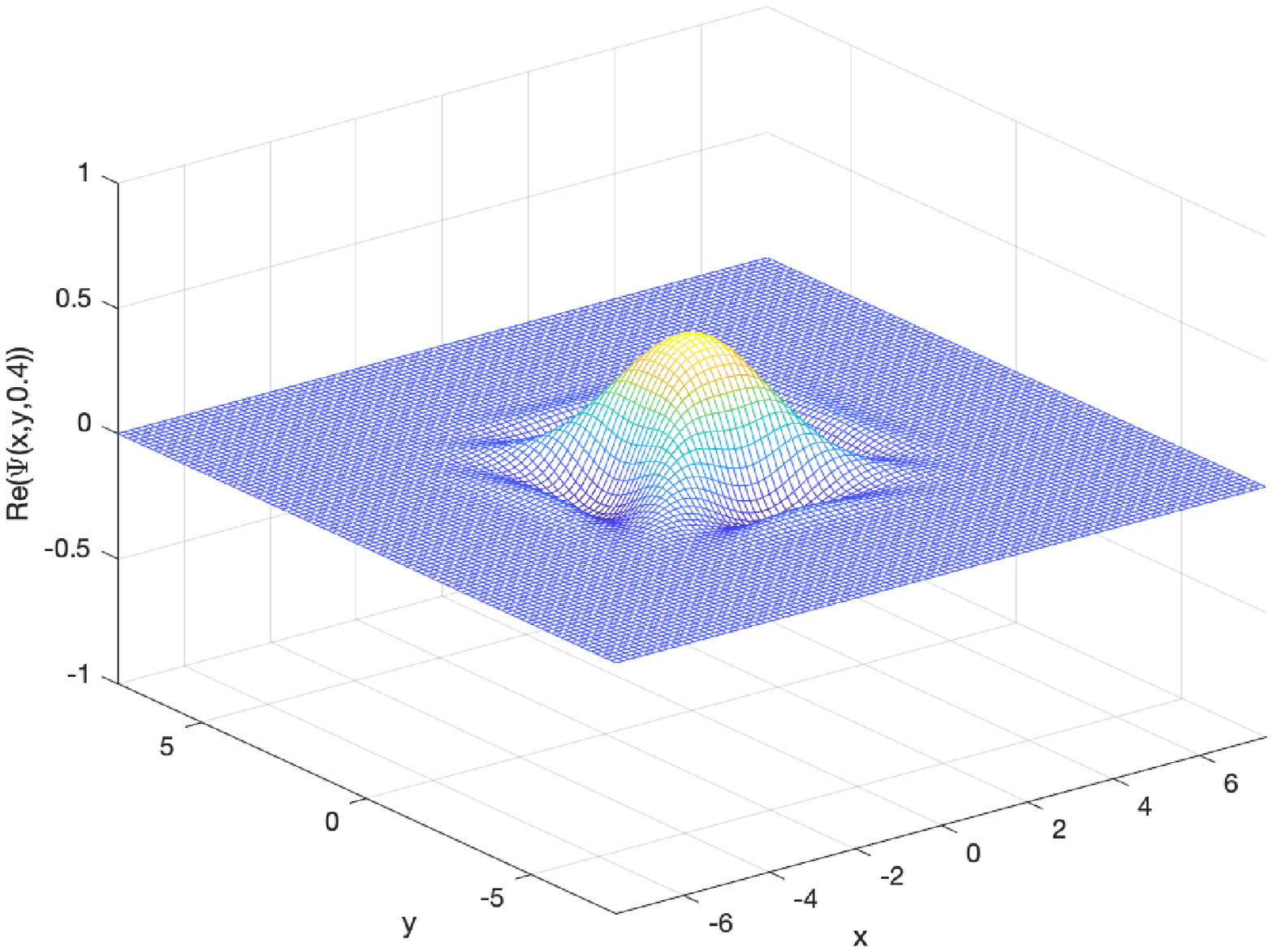}
\caption{Gaussian initial data (top) $\Psi(z,0)=e^{-|z|^{2}}$, and DS II evolution to time $t=0.4$ (bottom - $\mbox{Imag}(\Psi(z,t))$ on right, $\mbox{Re}(\Psi(z,t))$  on left) with $N=256$ Fourier modes, $L=4.56$, using the direct method describe in Section \ref{Sec:HybNum}.  }  
\label{fig:GaussianEvolution}
\end{figure}

Example 1:  The first example is the initial data $\Psi_{0} = e^{ - 
|z|^{2} } $.  We use the method from \cite{KM} to produce numerical 
approximations $\Psi(z,t)$ at time $t=0.4$. The solution is shown in 
Fig.~\ref{fig:GaussianEvolution}. It can be seen that the real and 
the imaginary part of the solution show oscillations. 

We study the convergence of the solution in dependence of the number 
$N$ of Fourier modes for $N = 2^{5}, 2^{6}, 2^{7}$, and $2^{8}$, on 
spatial domains $[-\pi L(N), \pi L(N)]$ that grow with $N$ 
(essentially $L(N)$ scales as $\sqrt{N}$).  For each reference 
solution for the same parameters, we compare to the solution to the DS II equation obtained using the direct solver as described in Section \ref{Sec:HybNum}.  The table below contains the relevant parameters, and maximum error between the two methods.
\begin{figure}[htb!]
\begin{tabular}{|c|c|c|c|}
\cline{1-4}
Fourier Modes & L & Error with classical method & Error with new method\\
\hline
$2^{5}$ & $1.65$ & \texttt{9.5e-05}& \texttt{1.2e-05} \\
\hline
$2^{6}$& $2.6$&1.5e-05 &\texttt{1.1e-07}\\
\hline
$2^{7}$ &$4.0$& 2.7e-06& \texttt{1.6e-10}\\
\hline
$2^{8}$& $6.15$&4.9e-07 &\texttt{1.2e-14} \\
\hline
\end{tabular}
\caption{Comparison of maximum error for the simulations with Gaussian initial data.  }  
\label{Fig:ParametersGauss}
\end{figure}

\begin{figure}[htb!]
\includegraphics[width=0.7\textwidth]{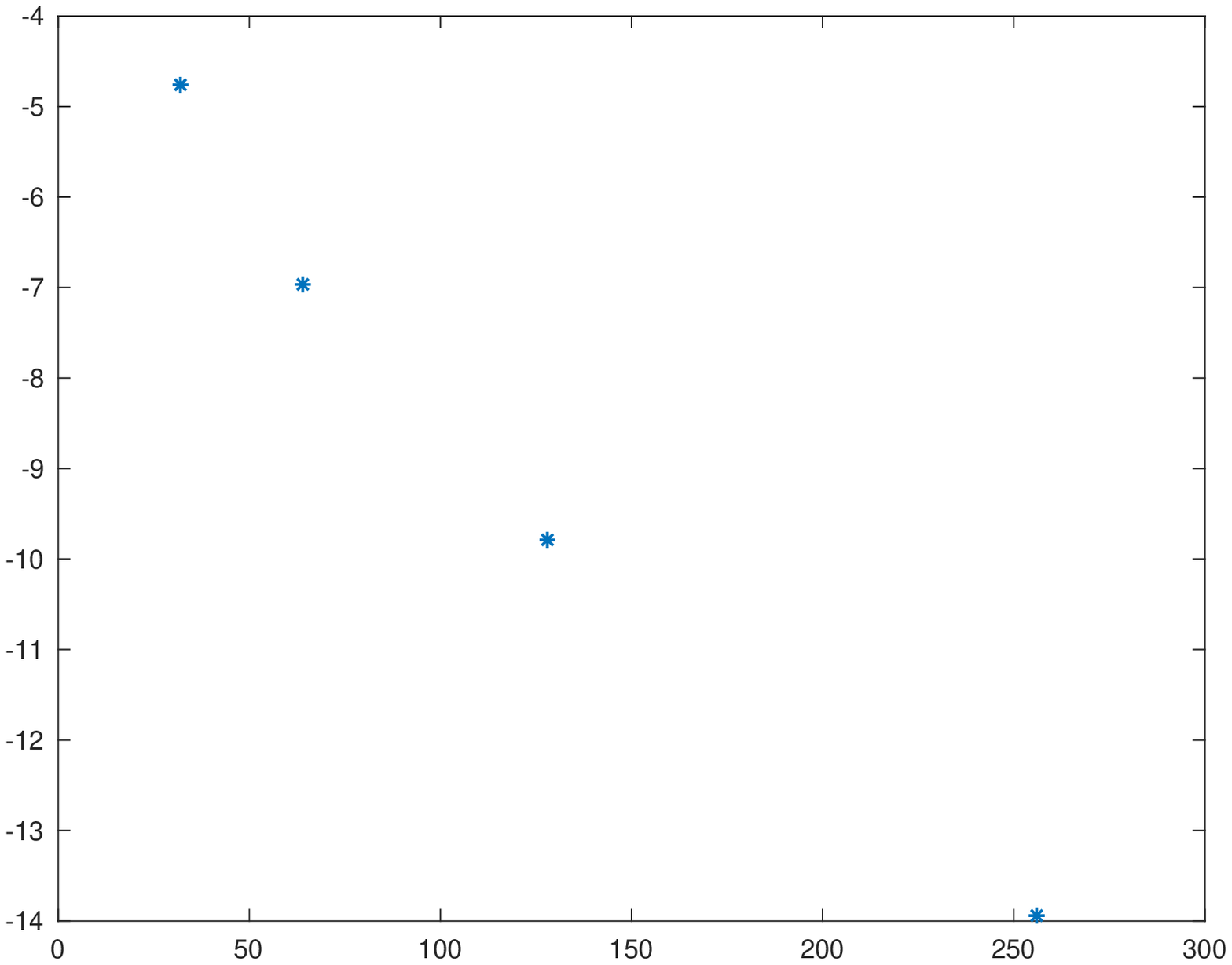}
\caption{Convergence comparison of the direct DS solver and the 
solution via inverse scattering, for $N=2^5, 2^6, 2^7$, and $2^8$, 
for defocussing DS II evolution of Gaussian initial data to time 
$t=0.4$ (see Fig.~\ref{fig:GaussianEvolution}).}  
\label{fig:GaussErrorPlot}
\end{figure}

The error of Table \ref{Fig:ParametersGauss} is plotted in 
Fig.~\ref{fig:GaussErrorPlot}. It can be seen that the error decays 
exponentially with $N$ as expected which means that spectral 
convergence is achieved. On the other hand the error for the 
classical codes for DS decreases only very slowly and essentially 
saturates at the order of $10^{-6}$. 

\begin{figure}[htb!]
\includegraphics[width=0.45\textwidth]{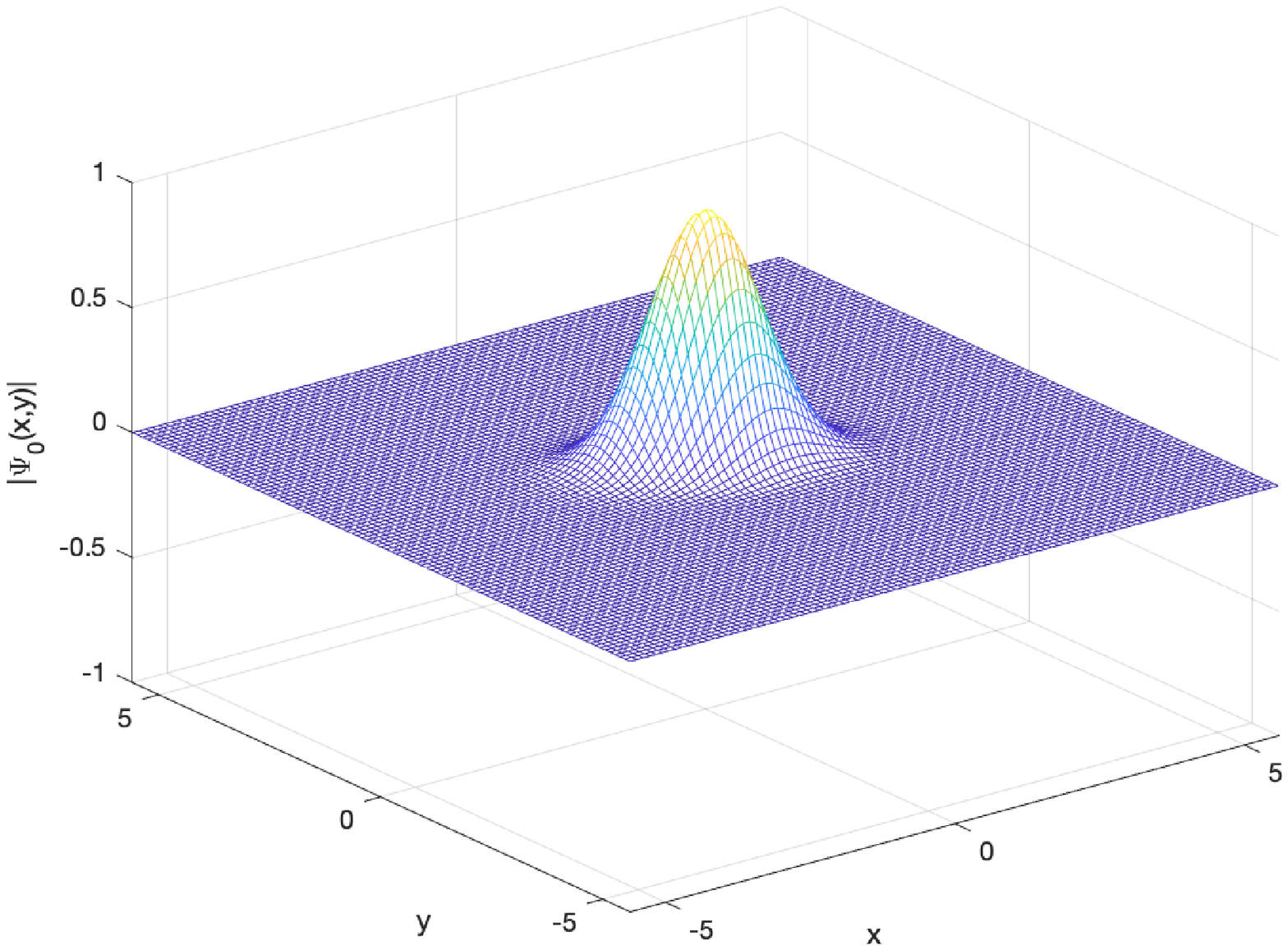}
\includegraphics[width=0.45\textwidth]{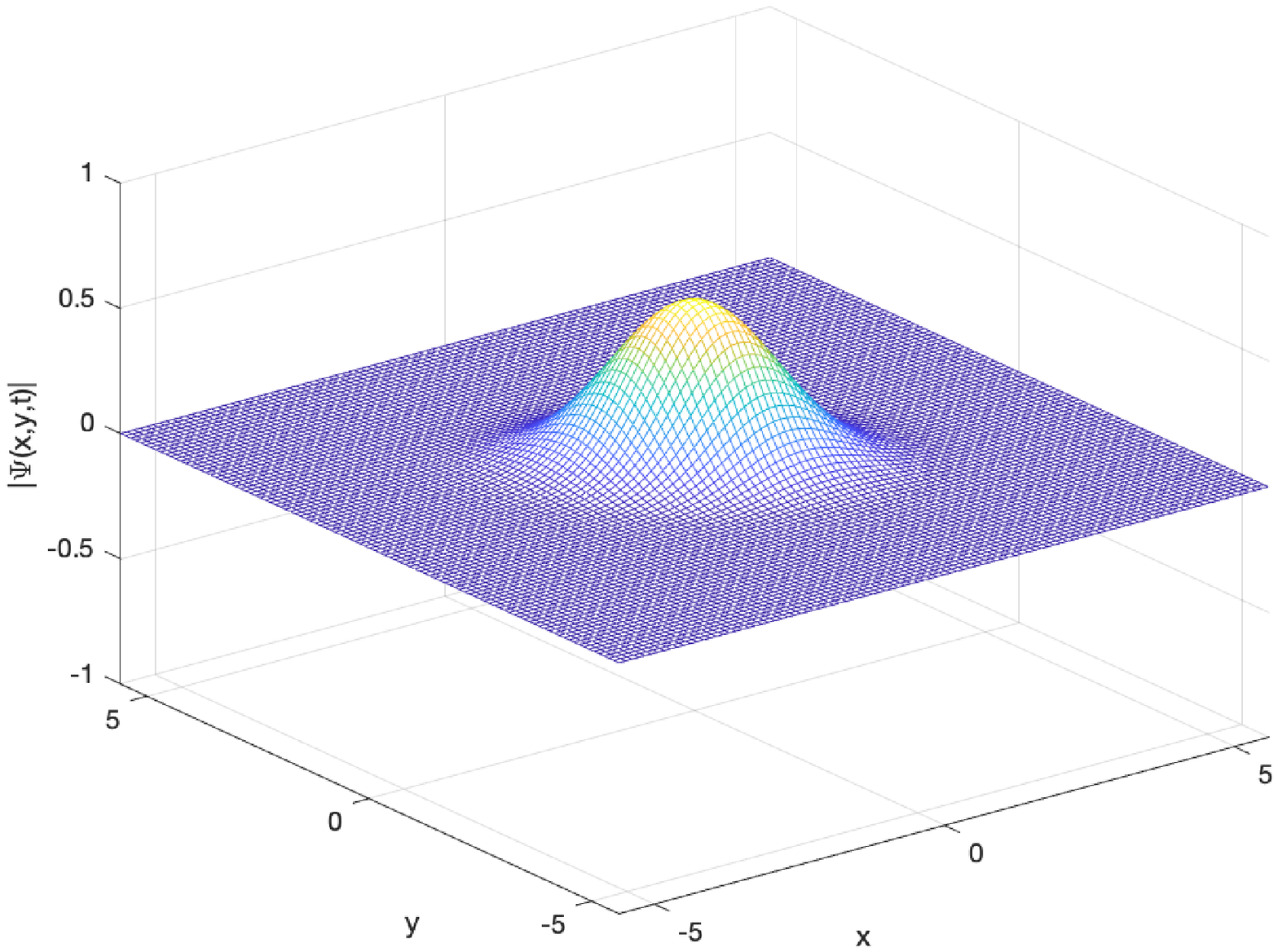} 
\includegraphics[width=0.45\textwidth]{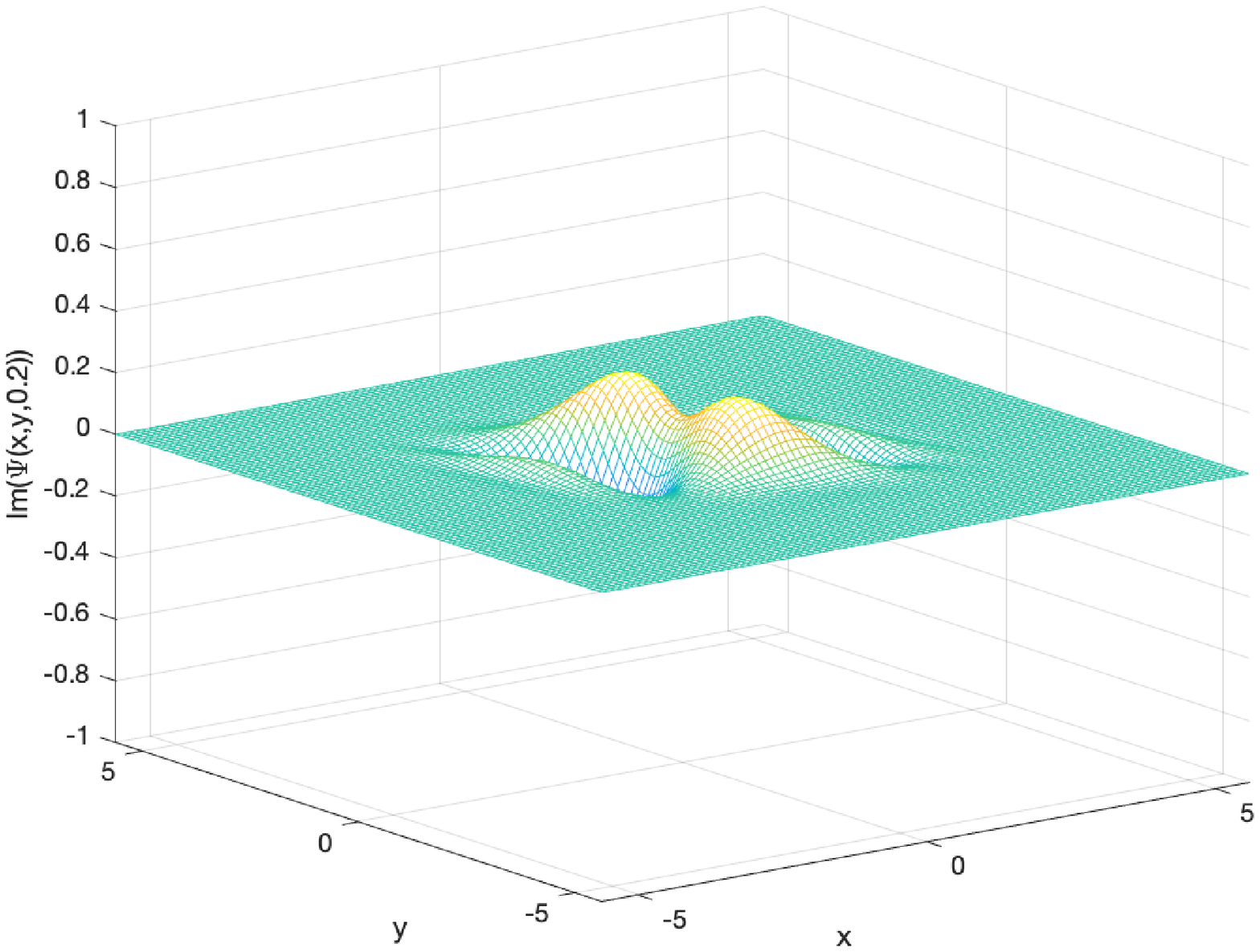}\includegraphics[width=0.45\textwidth]{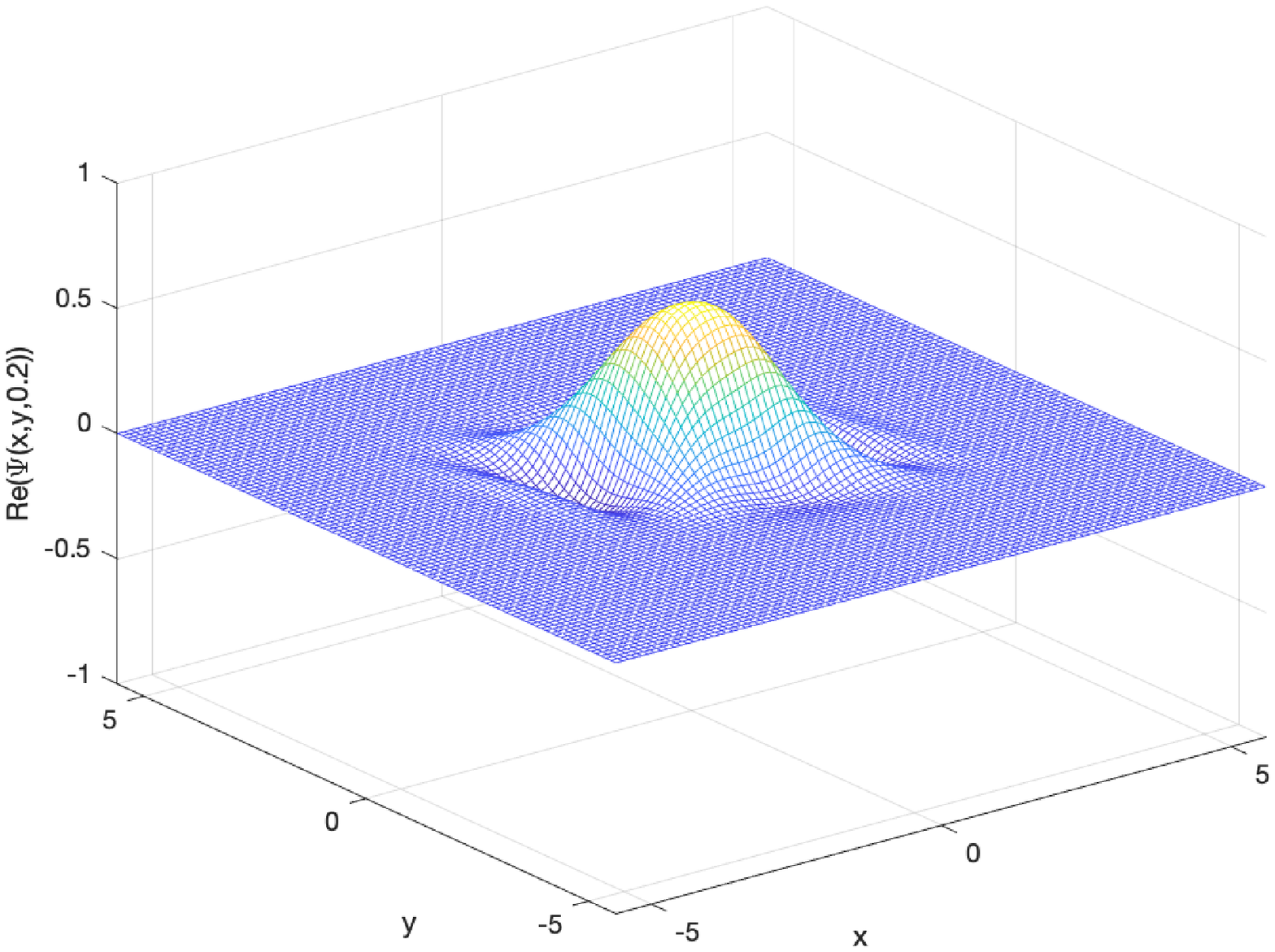}
\caption{Asymmetric initial data (top) $\Psi_{0} = 
e^{-\left( x^2 + xy + 2 y^2\right) }$, and DS II evolution to time 
$t=0.2$ (bottom - $\mbox{Imag}(\Psi(z,t))$ on right, $\mbox{Re}(\Psi(z,t))$  on left) with $N=256$ Fourier modes, $L=4.56$, using the direct method described in Section \ref{Sec:HybNum}.  }  
\label{fig:AsymmEvolution}
\end{figure}

Example 2:  The second example is the initial data $\Psi_{0} = 
e^{-\left( x^2 + xy + 2 y^2\right) }$, i.e., initial data in the 
Schwartz class without radial symmetry.  We again use the method from 
\cite{KM} to produce numerical approximations $\Psi(z,t)$, this time 
to time $t=0.2$. The solution is shown in 
Fig.~\ref{fig:AsymmEvolution}. 

The convergence of the numerical approach is again tested for $N = 
2^{5}, 2^{6}, 2^{7}$, and $2^{8}$, on spatial domains $[-\pi L(N), 
\pi L(N)]$ that scale with $\sqrt{N}$.  For each reference solution 
with the same parameters, we compare to the solution to the DS II equation obtained using the direct solver described in Section \ref{Sec:HybNum}.  The table below contains the relevant parameters, and maximum error between the two methods.
\begin{figure}[htb!]
\begin{tabular}{|c|c|c|c|}
\cline{1-4}
Fourier Modes & L & Error with classical method&Error with new method \\
\hline
$2^{5}$ & $1.41$&\texttt{7.6e-05}& \texttt{2.3e-05}\\
\hline
$2^{6}$& $2.12$&\texttt{1.3e-05}&\texttt{2.2e-07}\\
\hline
$2^{7}$ &$3.11$&\texttt{2.8e-06}& \texttt{7.1e-10} \\
\hline
$2^{8}$& $4.56$&\texttt{6.1e-07}&\texttt{3.2e-13} \\
\hline
\end{tabular}
\caption{Parameter values for the simulations with asymmetric initial data.  }  
\label{Fig:ParametersAsymm}
\end{figure}

\begin{figure}[htb!]
\includegraphics[width=0.7\textwidth]{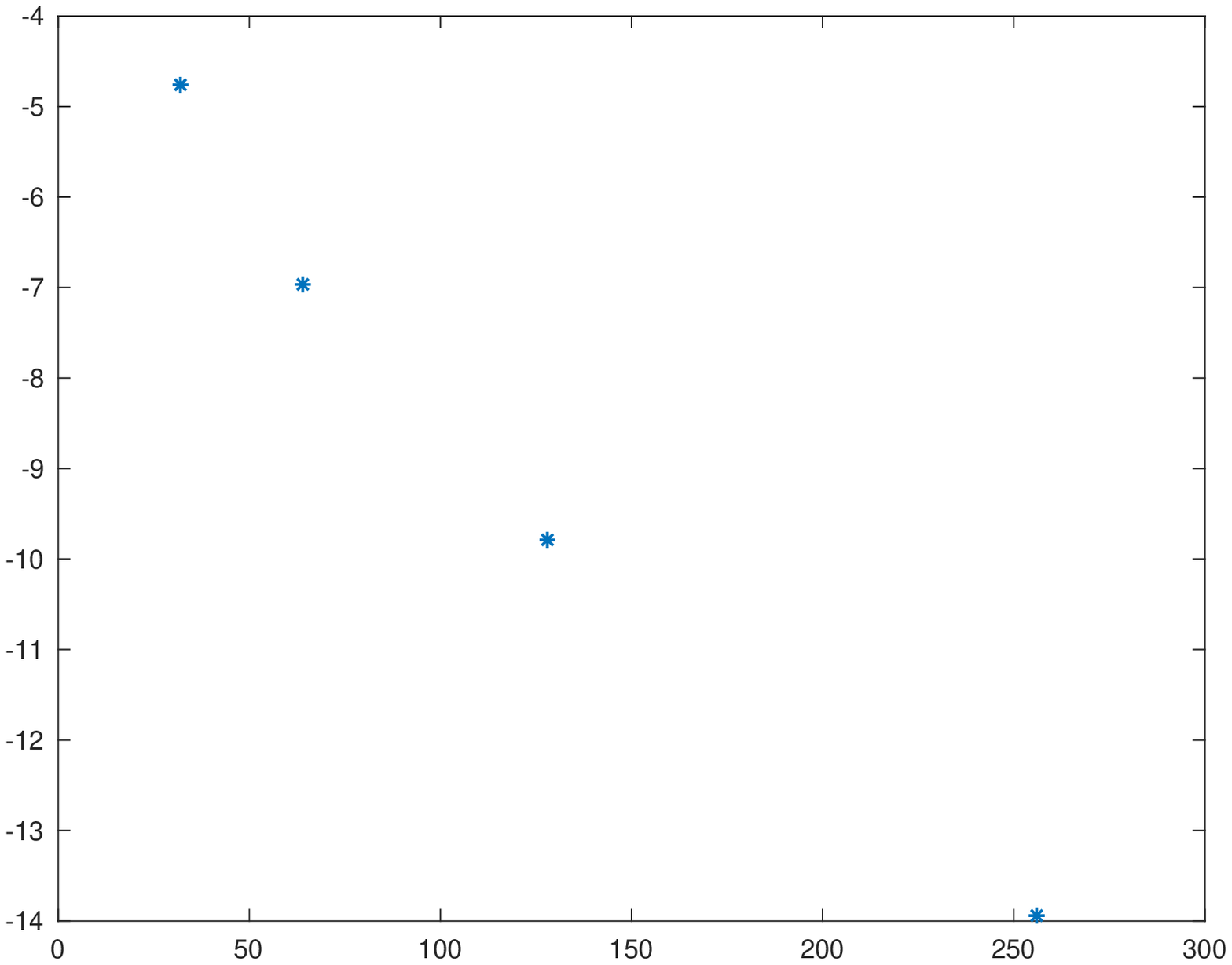}
\caption{Convergence comparison of two methods, $N=2^5, 2^6, 2^7$, 
and $2^8$, for defocussing DS II evolution of the asymmetric initial 
data to time $t=0.2$.  }  
\label{fig:AsymmErrorPlot}
\end{figure}

The exponential decay of the error is shown in 
Fig.~\ref{fig:AsymmErrorPlot}, whereas the error of the classical 
approach is once more essentially stagnant at the order of $10^{-6}$.

\section{Outlook}
\label{Sec:Conclusions}
In this paper we have presented a high precision numerical approach 
for DS II systems for Schwartz class initial data. Existing 
approaches based on a Fourier spectral method for the spatial 
dependence have been amended by regularizing the singular Fourier 
symbol of the nonlocal nonlinear term of the DS II systems via a hybrid method: 
the symbol is regularized via explicitly computable functions. Except 
for the numerical computation of these functions, which needs to be 
done only once for the whole time evolution, the same number of 
Fourier transforms (which are responsible for the main computational 
cost) are needed as in classical approaches. The accuracy of the 
method was tested via an independent numerical solution of the 
defocussing integrable DS II system by a numerical inverse scattering 
approach. It was shown that essentially machine precision can be 
reached with both methods. 

The Davey-Stewarson systems derived in the modulational regime of 
surface waves occur only in two-spatial dimensions. However the 
Zakharov-Rubenchik systems \cite{ZR} make physical sense also in 
three spatial dimensions. They may then degenerate to 
three-dimensonal Davey-Stewartson systems 
\begin{equation}\label{DS3D}
\begin{array}{ccc}
i\psi_{t} +\delta\psi_{xx}+\Delta_\perp 
\psi=(\chi|\psi|^2+\beta\phi)\psi,\\
\Delta\phi=\partial_{xx}|\psi|^2,
\end{array}
\end{equation}
where $\Delta_\perp=\partial_y^2+\partial_z^2$ and $\Delta 
=\partial_x^2+\partial_y^2+\partial_z^2$, and where $\delta$, $\chi$, 
and $\beta$ are real constants. 
This can be once more written as a nonlocal NLS equation, this time in 
3D, where the nonlocality corresponds to a singular Fourier symbol. 
The numerical study of these 3D DS II systems with a regularized 
symbol of the nonlocality along the lines of the present paper will 
be the subject of further research. 

\appendix
\section{Doubly periodic solution to DS II}\label{app}
\emph{by C. Kalla}\footnote{Universit\'e d'Orl\'eans, UFR Sciences,
MAPMO-UMR 6628, D\'epartement de Math\'ematiques,
Route de Chartres B.P. 6759 - 45067 Orl\'eans cedex 2, France,
E-mail address: caroline.kalla@univ-orleans.fr}

In this section we derive a a family of travelling wave solutions to the defocusing  DS II 
equation which are both periodic in $x$ and $y$. We only briefly 
touch the theory of algebro-geometric solutions to integrable 
equations. The reader is referred to \cite{BK} for an overview and a 
comprehensive list of references. 

It was shown in \cite{Kal} that the DS solutions in terms of 
multi-dimensional theta functions by Malanyuk \cite{Mal} can be 
written in the form
\begin{equation}
\Psi(x,y,t)=\sqrt{|q_{2}(a,b)|}\,\frac{\Theta(\z+\mathbf{r})}{\Theta(\z)}\,e^{ i(-N_{1}\Xi-N_{2}\eta+N_{3}\tfrac{t}{2})},
    \label{Psi}
\end{equation}
where
\[\Xi=x+iy, \qquad \eta=x-iy,\]
and
\[\z=i\mathbf{V}_{a}\,\Xi-i\mathbf{V}_{b}\,\eta+i(\mathbf{W}_{a}-\mathbf{W}_{b})\,\tfrac{t}{2};\]
here the multi-dimensional theta function is defined as the series
\begin{equation}
    \Theta(\z) = \sum_{\mathbf{n}\in\mathbb{Z}^{g}}^{}\exp(\langle 
    \mathbf{n}|\mathbb{B}\mathbf{n}\rangle/2+\langle\mathbf{n}|\mathbf{z}\rangle)
    \label{theta},
\end{equation}
which is uniformly convergent since $\mathbb{B}$ is a Riemann matrix, 
i.e., satisfies $\mathbb{B}^{t}=\mathbb{B}$ and has a negative 
definite real part; $\langle \cdot|\cdot\rangle$ defines the 
euclidean scalar product, $a$, $b$ are points on a Riemann surface of 
genus $g$, and $q_{2}$, $N_{1}$, $N_{2}$, $N_{3}$, $\mathbf{r}$, $\mathbf{V}$, 
$\mathbf{W}$ are certain integrals respectively their periods on this 
Riemann surface (see \cite{Kal}). We do not specify them here since 
they will be just parameters to be chosen in a way to get 
periodic solutions. 

We now try to identify some of the simplest periodic solutions 
amongst (\ref{Psi}). 
Let $L_{x}$ and $L_{y}$ be the $x$-period and $y$-period respectively and let
\[\tilde{\Psi}(x,y,t)=\Psi(x+k_{1}L_x,y+k_{2}L_y,t),\quad 
k_{1},k_{2}\in \mathbb{Z}.\]
We want solutions such that
\[\tilde{\Psi}(x,y,t)=\Psi(x,y,t)\quad  \forall x,y,t,k_{1},k_{2}.\] 
We can see that
\[\tilde{\Psi}(x,y,t)=\sqrt{|q_{2}(a,b)|}\,\frac{\Theta(\z+\mathbf{r}+\mathbf{L_{1}})}{\Theta(\z+\mathbf{L_{1}})}\,e^{i(-N_{1}\Xi-N_{2}\eta+N_{3}\tfrac{t}{2})}\,e^{L_{2}},\]
where
\[\mathbf{L_{1}}=k_{1}L_{x}(\mathbf{V}_{a}-\mathbf{V}_{b})i-k_{2}L_{y}(\mathbf{V}_{a}+\mathbf{V}_{b}),\]
\[L_{2}=-k_{1}L_{x}(N_{1}+N_{2})i+k_{2}L_{y}(N_{1}-N_{2}).\]
\\
In order to obtain periodic solutions, we need that
\begin{equation}
\mathbf{L_{1}}=2i\pi\mathbf{n}+\B\mathbf{m}, \qquad \mathbf{n},\mathbf{m}\in\Z^{2}, \label{1}
\end{equation}
which implies
\[\tilde{\Psi}(x,y,t)=\sqrt{|q_{2}(a,b)|}\,\frac{\Theta(\z+\mathbf{r})}{\Theta(\z)}\,e^{i(-N_{1}\Xi-N_{2}\eta+N_{3}\tfrac{t}{2})}\,e^{L_{2}-\left\langle \mr,\mathbf{m}\right\rangle}.\]
Hence, if (\ref{1}) is satisfied, we need that
\begin{equation}
\Re(L_{2})=\left\langle \Re(\mr),\mathbf{m}\right\rangle. \label{2}
\end{equation}
\\
Now, recall that in order to get a regular solution to the defocusing 
DS II equation, one must have
\begin{equation}
\mathbf{V}_{b}=-\overline{\mathbf{V}_{a}}, \qquad N_{2}=\overline{N_{1}}. \label{3}
\end{equation}
This implies 
\begin{equation}
\mathbf{L_{1}}=2i\left(k_{1}L_{x}\,\Re(\mathbf{V}_{a})-k_{2}L_{y}\,\Im(\mathbf{V}_{a})\right),\label{4}
\end{equation}
and
\begin{equation}
L_{2}=2i\left(-k_{1}L_{x}\,\Re(N_{1})+k_{2}L_{y}\,\Im(N_{1})\right).\label{5}
\end{equation}
Now, if we assume that (for defocusing DS II)
\[ \overline{\B}=\B,\]
we deduce from (\ref{1}) and (\ref{4}) that
\[\mathbf{m}=0,\]
and (\ref{2}) is verified. Hence, always by (\ref{1}) and (\ref{4}), we search $L_{x},\,L_{y}\in\R$ such that for a given $\mathbf{n}\in\Z^{2}$
\begin{equation}
k_{1}L_{x}\,\Re(\mathbf{V}_{a})-k_{2}L_{y}\,\Im(\mathbf{V}_{a})=\pi \mathbf{n}. \label{6}
\end{equation}
Since this equation has to hold for all $k_{1},k_{2}\in\mathbb{Z}$, 
we get
\begin{align}
    n_{1}\Re(V_{a,1}) &= n_{2}\Re(V_{a,2}),
    \nonumber\\
     m_{1}\Im(V_{a,1}) & =m_{2}\Im(V_{a,2})
    \label{cond1}
\end{align}
for some $m_{1},m_{2},n_{1},n_{2}\in \mathbb{Z}$. A possible choice 
is  $\Re(V_{a,1})=\Im(V_{a,2})=0$. If $ \Re(V_{a,1})$, $\Im(V_{a,1})$ do 
not vanish, one has
$$L_{x}=\frac{\pi n_{1}}{\Re(V_{a,1})},\quad L_{y}=-\frac{\pi 
m_{1}}{\Im(V_{a,1})}.$$

The remaining condition to be satisfied is (\ref{5}) which reads
\begin{equation}
    k_{2}L_{y}\Im (N_{1})-k_{1}L_{x}\Re(N_{1})=k\pi, \quad k\in\mathbb{Z}
    \label{5a}.
\end{equation}
Since this has to hold for all $k_{1}$, $k_{2}$, this is again 
equivalent to two equations,
\begin{equation}
    \frac{\Re(N_{1})}{\Re(V_{a,1})}=l_{1},\quad 
    \frac{\Im(N_{1})}{\Im(V_{a,1})}=l_{2}, \quad 
    l_{1},l_{2}\in \mathbb{Z}.
    \label{5b}
\end{equation}

The task is thus to find a point $a$ on an M-curve of genus 2 (which 
means a hyperelliptic curve with only real branch points) such that 
(\ref{cond1}) and (\ref{5b}) are satisfied and such that the solution is 
non trivial. Thus we consider the hyperelliptic curve 
\begin{equation}
 \mu^{2}=\prod_{n=1}^{6}(\lambda-\lambda_{n}),\quad \lambda_{n}\in 
\mathbb{R}  . 
    \label{bp}
\end{equation}
The vector $\mathbf{V}$ on this curve has the form $\mathcal{A}(1,a)^{t}/\mu(a)$, 
where $\mathcal{A}$ is the inverse of the matrix of $a$-periods of 
the holomorphic one-forms times 
$2\pi i$. One possibility is to choose the curve and the point $a$ such 
that $V_{a,1}=\bar{V}_{a,2}$. In this case equation (\ref{cond1}) is obviously 
satisfied. This can be achieved on a curves with a symmetry with 
respect to $\lambda\mapsto -\lambda$ and a point $a$ projecting to the imaginary 
axis. But this leaves two conditions (\ref{5b}) to be satisfied with 
a single real parameter which implies that though a whole family of 
such solutions exist, there must be a special relation between some 
of its parameters. 

In order to construct a particular such solution, we choose the 
branch  points in (\ref{bp}) to be $\alpha,2,1,-1,-2,-\alpha$ with 
$\alpha>2$ and $\mathbf{V}_{a}=\beta i[1,-1]$ iteratively in 
order to satisfy conditions (\ref{cond1}) for given $n_{1}/n_{2}$, 
$m_{1}/m_{2}$. This is done with some initial guess with the 
optimization algorithm \cite{lagarias} implemented in Matlab as the command 
\emph{fminsearch}. The quantities on the hyperelliptic curve are 
computed with the code \cite{FK}. To give a concrete example, we 
choose $n_{2}/n_{1}=1$ and $m_{2}/m_{1}=-4$ and get 
$\alpha=2.870255599870804$, $\beta=2.119032837086884$ and thus the 
solution shown in Fig.~\ref{figexact}. Internal tests of the code 
\cite{FK} indicate that the solution satisfies DS II to the order of 
$10^{-14}$. 

\section{Integrability of the DS II equation}

In this section of the Appendix, we provide a brief summary of the scattering and inverse scattering theory that provides a solution procedure for the integrable case of the defocusing DS II equation.

Here we follow the exposition in \cite{KM}.  From the initial data 
$\Psi_{0}(x)=\Psi(x,0)$, we determine the \emph{reflection 
coefficient} (the scattering transform of the initial data, sometimes 
also called a nonlinear Fourier transform) by solving the following system of  linear elliptic partial differential equations
\begin{eqnarray}
\label{eq:specprob2}
\pmtwo{\dbar}{0}{0}{\partial} \psi = 
\frac{1}{2}\pmtwo{0}{\Psi_{0}(z)}{\overline{\Psi_{0}(z)}}{0} \psi .
\end{eqnarray}
The operators $\partial $ and $\overline{\partial}$ are defined via
\begin{eqnarray*}
\partial = \frac{1}{2} \left(
\frac{\partial}{\partial x} - i \frac{\partial}{\partial y}
\right), \ \ \ \ \ 
\overline{\partial} = \frac{1}{2} \left(
\frac{\partial}{\partial x} + i \frac{\partial}{\partial y}
\right).
\end{eqnarray*}

We seek a column vector 
$\psi = \psi(z,k) = \left( \begin{array}{c}
\psi_{1} \\
\psi_{2} \\ \end{array}
\right)$ solving (\ref{eq:specprob2}), with the following asymptotic behavior as $|z| \to \infty$:
\begin{eqnarray*}
&&\lim_{|z| \to \infty} \psi_{1} e^{-kz/\epsilon} = 1, \\
&& \lim_{|z| \to \infty} \psi_{2} e^{-\overline{k} \overline{z} / \epsilon} = 0 \ .
\end{eqnarray*}
Here the quantity $k$ is a complex parameter, $k=k_{1} + i k_{2} $ with $(k_{1},k_{2})\in\mathbb{R}^2$, which plays the role of a spectral variable.
In keeping with the literature on inverse problems, we refer to the quantity $\psi$ as a \emph{complex geometric 
optics} (CGO) solution.   

Given initial data $\Psi_{0}(z)$ for the DSII equation, the CGO 
solution $\psi$ is uniquely determined, and the reflection coefficient $r=r^{\epsilon}(k)$, is obtained from the sub-leading term in the asymptotic expansion of $\psi$ as $z \to \infty$, via
\begin{eqnarray}\label{reflec}
\psi_{2} e^{ - \overline{k} \overline{z}} = \frac{\bar{r^{\epsilon}}(k)}{2\overline{z}} + \mathcal{O} \left( \frac{1}{|z|^{2}} \right) \ .
\end{eqnarray}
The mapping from $\Psi_0$ to $r$ is a transformation from the potential 
$\Psi_{0}(x,y)$ (again a function of two real variables) to a function 
$r(k_{1}, k_{2})$, which extends to a Lipschitz continuous and 
invertible mapping on the function space $L^{2}\left( \mathbb{C} 
\right)$ (see \cite{Perry2012,NRT}, and the references contained therein).  

Now, if $\Psi =\Psi(x,y,t)$ evolves according to the DSII equation (\ref{DSgen}), then the reflection coefficient evolves according to
\begin{eqnarray*}
r=r(k,t) = r(k,0) e^{ \frac{-it}{4 }  \left( k^{2} + \overline{k}^{2} \right) } \ .
\end{eqnarray*}
It is well-known that the inverse problem, reconstructing the potential $q(x,y,t,\epsilon)$ from the reflection coefficient $r(k,t)$, is also a D-bar problem, this time with respect to the complex variable $k$.  Indeed, setting 
\begin{equation}
\phi_1:=
\ee^{-kz}\psi_1\quad\text{and}\quad
\phi_2:=
\ee^{-\overline{k} \overline{z}}\psi_2 \ , 
\label{Phi}
\end{equation}
one may verify that for each  $z \in \mathbb{C}$,
\begin{equation}
\overline{\partial}_k \phi_1=\tfrac{1}{2}e^{(\overline{k} 
\overline{z} - k z )} \overline{r(k,t)}\phi_2, \ \ \ \ \ \ \ \ 
\partial_k \phi_2 =\tfrac{1}{2}e^{-(\overline{k} \overline{z} - k z )}r(k,t)\phi_1
\label{eq:dbar-k}
\end{equation}
where,
\begin{equation*}
\overline{\partial}_k:=\frac{1}{2}\left(\frac{\partial}{\partial k_{1}}+\ii\frac{\partial}{\partial k_{2}}\right), \ \ \ \ \ \partial_{k}:=\frac{1}{2}\left(\frac{\partial}{\partial k_{1}}-\ii\frac{\partial}{\partial k_{2}}\right),
\end{equation*}
and the quantities $\phi_{1}$ and $\phi_{2}$ possess the following asymptotic behavior:
\begin{equation}
\lim_{|k|\to\infty}\phi_1(k;z,t)=1\quad\text{and}\quad
\lim_{|k|\to\infty}\phi_2(k;z,t)=0.
\label{eq:Phi-asymp}
\end{equation}
The functions $\phi_{1}$ and $\phi_{2}$ are uniquely determined by 
this elliptic PDE system and boundary conditions.  The potential $q(x,y,t,\epsilon)$ is similarly determined through the asymptotic behavior as $|k| \to \infty$:
\begin{eqnarray*}
\phi_{2} =  \frac{\overline{q(x,y,t,\epsilon)}}{2 k } + \mathcal{O} \left( |k|^{-2}\right) \ .
\end{eqnarray*}

%


\end{document}